\documentclass[preprint,12pt,3p]{elsarticle}

\usepackage{amssymb}
\usepackage{gensymb}
\usepackage{graphicx}
\graphicspath{{Images/}}

\journal{arXiv}

\begin{document}

\begin{frontmatter}

\title{Properties of a Two Dimensional Model of RNA Folding}
%\tnotetext[label0]{This is only an example}

\author[label1]{Ben Y. Maron}
\address[label1]{Williams College, 880 Main St, Williamstown, MA 01267, USA}

\ead{benmaron@gmail.com}

\begin{abstract}
Ribonucleic Acid (RNA) can fold into shapes that perform functions in the cell. These foldings are governed by Watson-Crick base pairing, namely Adenine to Uracil and Cytosine to Guanine (A-U and G-C). The properties of the H-P (hydrophobic-hydrophilic) model of protein folding has been well studied in the two dimensional orthogonal case, and we attempt to achieve similar results. We prove that (1) there is an infinite family of even-length sequences with unique optimal foldings, (2) there are two ideal foldings for an even length sequence and given a sequence is is quickly verifiable if both, one, or neither are optimal, (3) finding an optimal foldings under this model is NP-hard, and (4) we give a constant-factor approximation algorithm for optimally folding RNA sequences.
\end{abstract}


\end{frontmatter}

% \linenumbers

\section{Introduction}
\label{sec1}

RNA and DNA are each composed of four bases which bond in pairs. DNA uses Adenine, Thymine, Guanine, and Cytosine (A, T, G, and C) to encode information, and RNA replaces Thymine with Uracil (U). They bond with each other in the pattern G-C and A-T/U \cite{rich}. However, while DNA typically folds into a stable double helix, RNA can fold into a wide variety of shapes depending on its sequence. This allows RNA to function as more than information storage, and take an active role in the cell by catalyzing reactions and functioning as an enzyme \cite{rnazymes}. One of the most fundamental ribozyme is the ribosome, which carries out protein synthesis \cite{rrna}.

While RNA folding has remained relatively unexplored from a mathematical perspective, significant results have been proven for protein folding, specifically the hydrophobic-hydrophilic (H-P) model, which assumes a protein is a chain where each node is either hydrophobic or hydrophilic. An energy function is defined by counting the number of adjacent hydrophobic nodes which are not adjacent in the chain. An optimal folding is one which maximizes this number for a given sequence \cite{hpmodel}. Assuming the chain is orthogonally embedded in two dimensions, is has been shown that there is a sequence of any doubly even length with a unique optimal folding \cite{uniqueness}, and that in general finding an optimal folding is NP-hard \cite{hardness}. Furthermore, a constant-factor approximation has been found which guarantees a folding within one-third of optimal \cite{approximation}. This paper aims to prove similar results for a related model of RNA folding. As a result, much of the terminology and setup will be shared with papers on the H-P model.

\section{The Watson-Crick Model}
\label{sec2}

RNA is a polymer whose monomers are nucleotides. There are four nucleotides, A, U, G, and C, and generally they bond A-U and G-C. Occasionally non-Watson-Crick base pairing occurs \cite{nonwatsoncrick}, or nucleotides are modified in a way which may change their properties \cite{modification}. For the purposes of maintaining a simple model, we assume that these do not occur. 

A chain refers to a linkage representing a sequence of nucleotides with each node assigned a base, and a folding is a non-self intersecting embedding of a chain in a two dimensional square lattice. In all figures, the chain is denoted by a black line and nucleotides are shown as circles on the nodes of the chain. Their assignment is as follows: G-black, C-white, A-horizontal line, U-vertical line. Bonds can form wherever two complementary bases are adjacent in the folding but not in the chain, and are denoted by dashed lines. Each node can only participate in one bond because in nature each nucleotide can only bond to one other. An optimal folding has the maximum number of bonds for any folding of a given chain, and a uniquely optimal folding is one where no other folding (disregarding symmetries) attains the same number of bonds.

\section{Uniquely Foldable Chains}
\label{sec3}

\begin{figure}[h]
\includegraphics[scale=0.15]{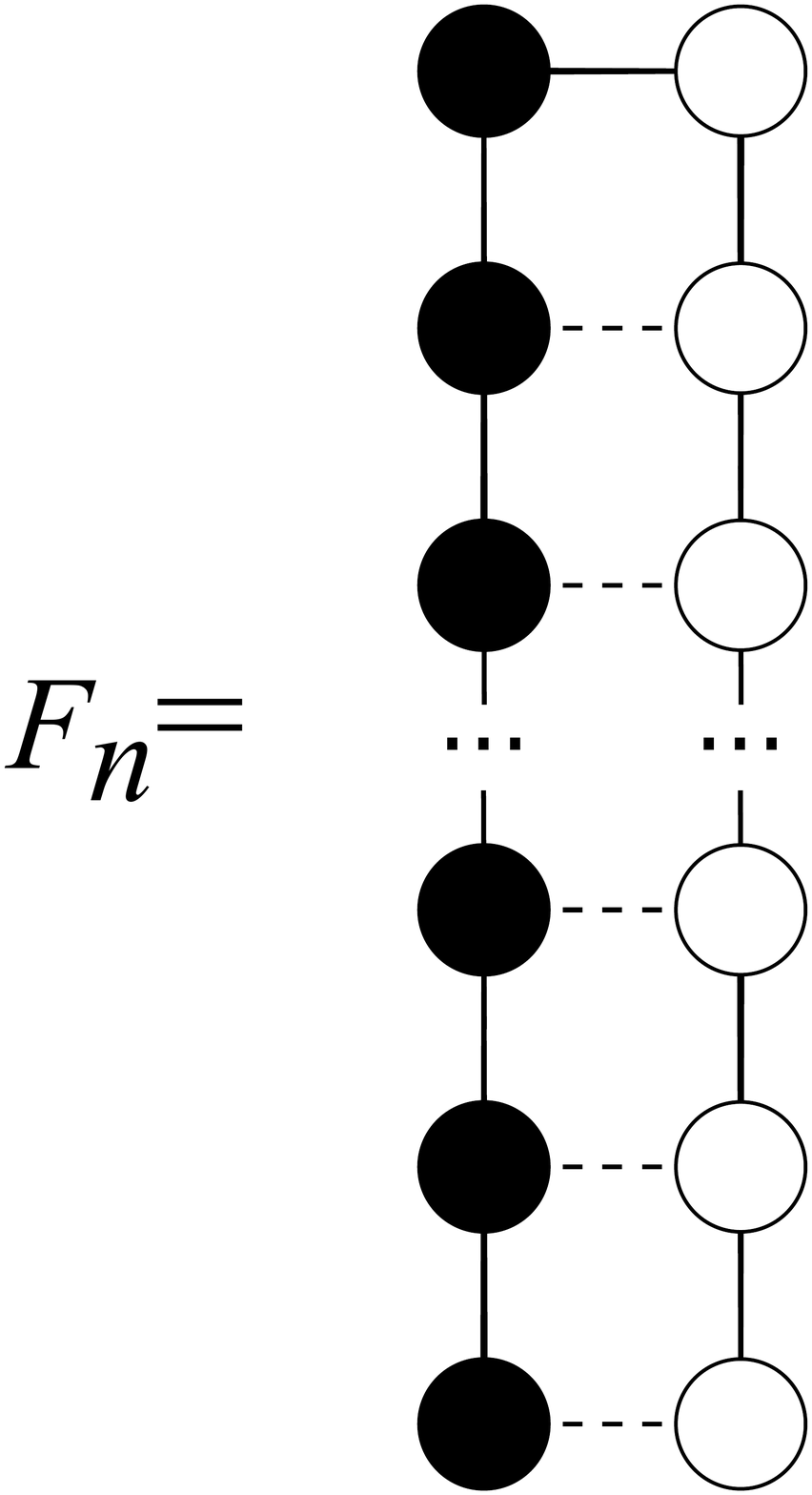}
\centering
\caption{Illustration of $F_n$.}
\label{fn}
\centering
\end{figure}

\noindent
\textbf{Theorem 1.} \textit{For $n>3$, the chain $S_n=G^nC^n$ has a unique optimal folding $F_n$ (Figure \ref{fn})}

\bigskip
\noindent
\textbf{Lemma 2.} \textit{The upper bound on the number of bonds in a folding of any chain with length $2n$ is $n-1$.}

\bigskip
\noindent
\textbf{Proof.} Consider the bounding box of the folding. Pick one edge and examine one of its extreme nodes. If it is an endpoint, it may form a bond with a medial node.

\begin{figure}[h]
\includegraphics[scale=0.13, angle=90]{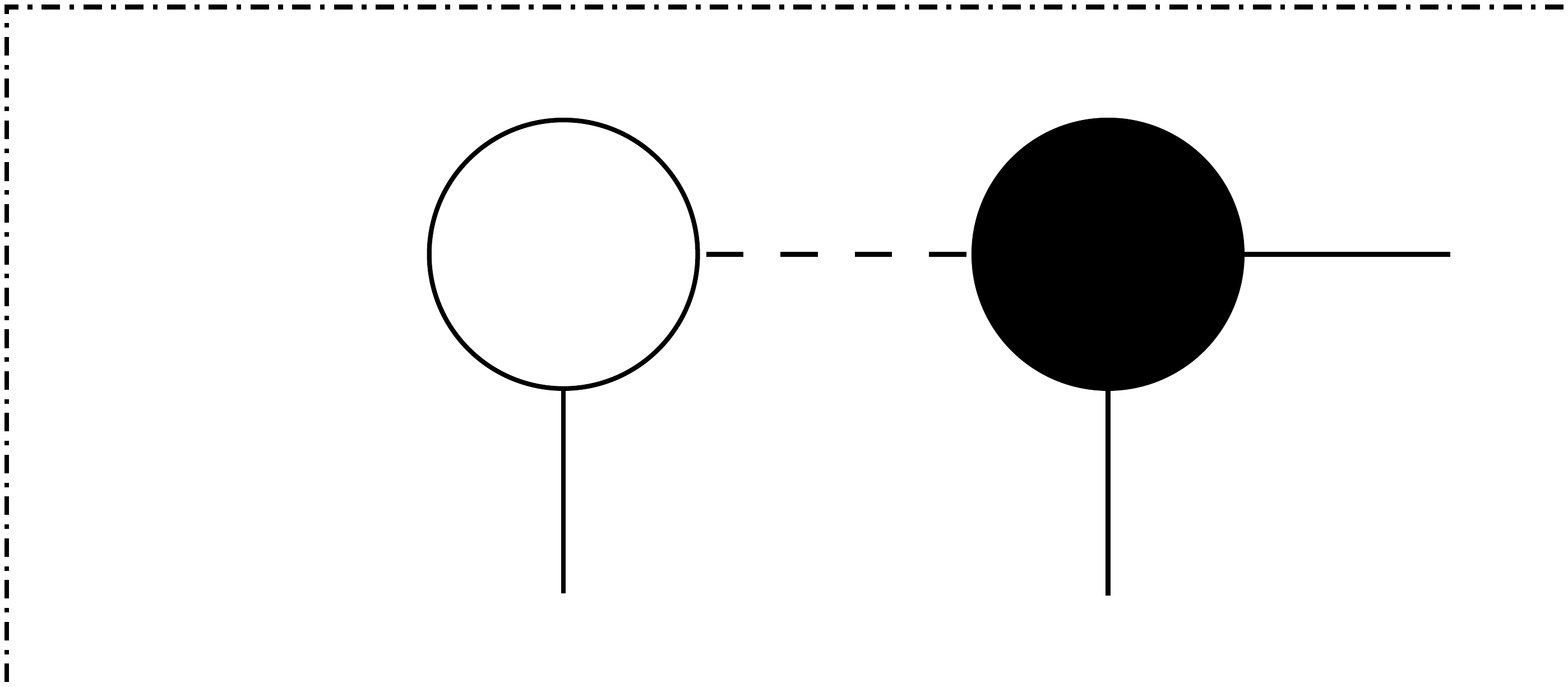}
\centering
\caption{An endpoint can form a bond as an extreme node (bounding box shown as dot-dashed line).}
\centering
\end{figure}

If it is not an endpoint, it cannot form any bonds because the only two adjacent locations are filled by nodes that are directly adjacent to it in the chain.

\begin{figure}[h]
\includegraphics[scale=0.13, angle=90]{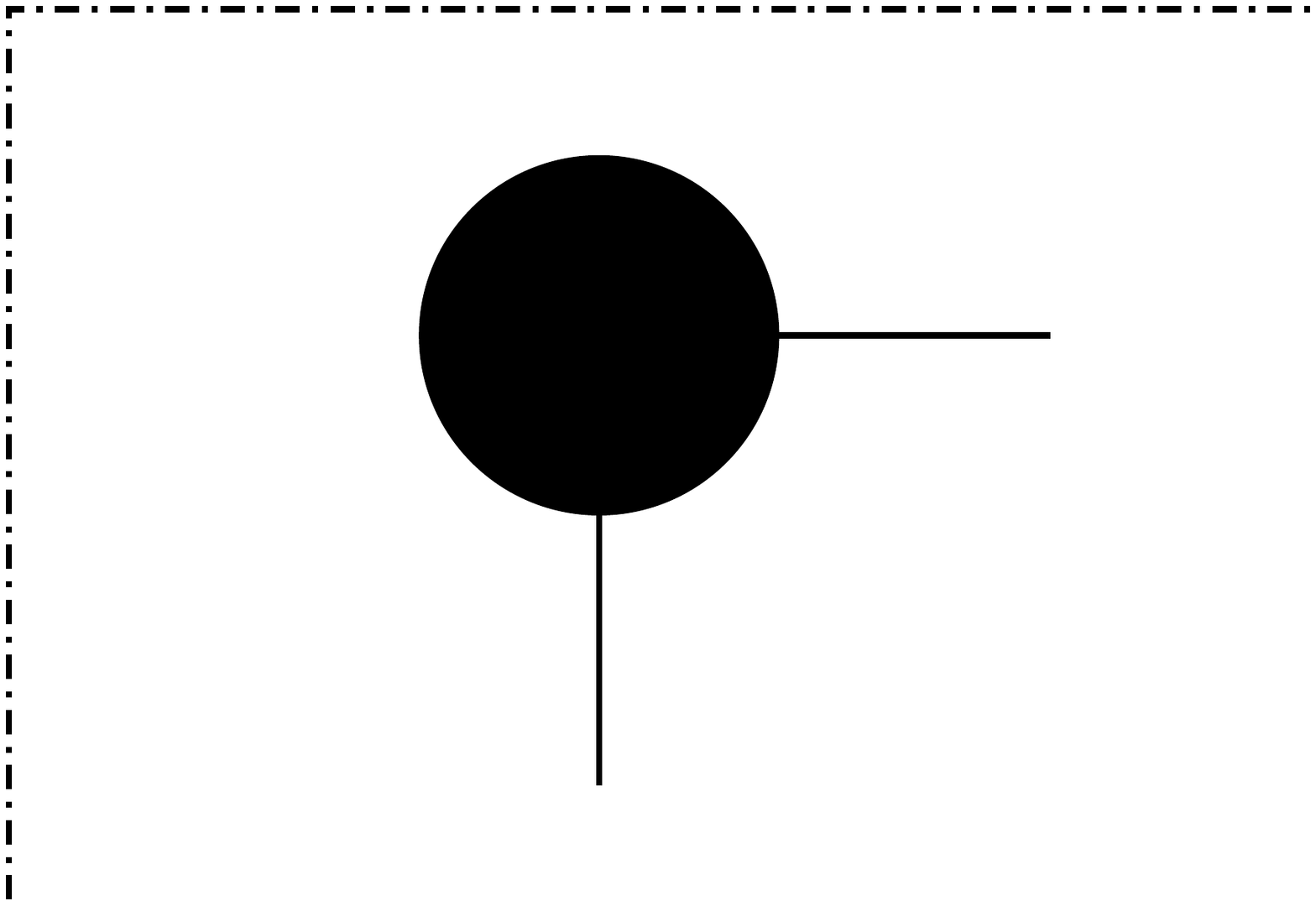}
\centering
\caption{A non-endpoint can not form a bond as an extreme node.}
\centering
\end{figure}

There are eight extreme nodes, two on each side of the bounding box, but that number can be reduced. If a node is in the corner of the bounding box, it can be the extreme node for two faces. Additionally, if there is only one node on a side, it is both extremes for that side. However, the only nodes that could fit that description are the endpoints, and (as will be shown) the endpoints need to bond in order to achieve the upper bound. Therefore, placing one node in each corner of the bounding box minimizes the number of unbound nodes, and making two of the corner nodes endpoints means that only two nodes must be unbound. If there are $2n$ total nodes and $2$ are unbound, there are $2n-2$ nodes that can bond, which would form a maximum of $n-1$ bonds. $\Box$

\bigskip
\noindent
\textbf{Fact 3.} \textit{$F_n$ has $n-1$ bonds.}

\bigskip
\noindent
\textbf{Fact 4.} \textit{The two unbound nodes of $F_n$ must be one C and one G.}

\bigskip
\noindent
\textbf{Proof.} If the nodes were the same base, there would not be an equal number of remaining C and G nodes so they would not all be able to bond. $\Box$

\bigskip
\noindent
\textbf{Fact 5.} \textit{The two unbound nodes of $F_n$ cannot be the endpoints.}

\bigskip
\noindent
\textbf{Proof.} Lemma 2 uses the fact that the endpoints can fill the corners of the bounding box but also form bonds in order to reach the upper bound. $\Box$

\bigskip
\noindent
\textbf{Lemma 6.} \textit{Diagonally opposite corner nodes can not be the same base.}

\bigskip
\noindent
\textbf{Proof.} Assume the contrary. The endpoints must be in two adjacent corners because otherwise the assumption would not hold, and the other two corner nodes must have opposite assignments by Fact 3. Starting from the G endpoint, there must be a chain of solid G's that connect the two opposite G corners before the chain reaches the other two because of the assignment of $S_n$. However, there is no way for the chain to then reach the other two corners without self-intersection. Therefore the assumption is false and diagonally opposite corner nodes must have opposite assignments. $\Box$

\begin{figure}[h]
\includegraphics[scale=0.13, angle=90]{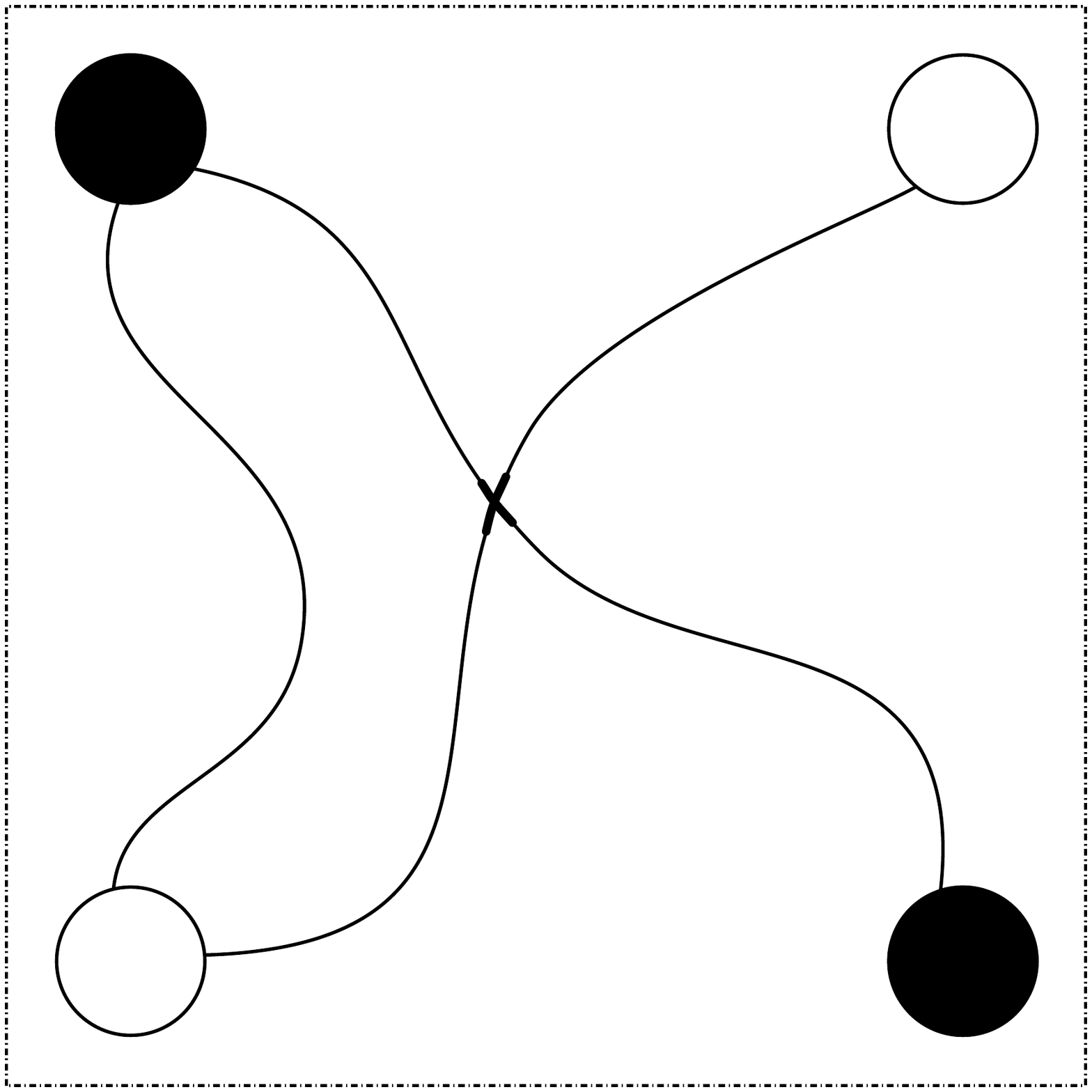}
\centering
\caption{Illustrating the proof of Lemma 6.}
\centering
\end{figure}

\bigskip
\noindent
\textbf{Lemma 7.} \textit{If $n>3$, the two endpoints must share a face.}

\bigskip
\noindent
\textbf{Proof.} Assume the contrary. Due to the assumption and Lemmas 1 and 6, the optimal folding must look like Figure 5.

\begin{figure}[h]
\includegraphics[scale=0.13]{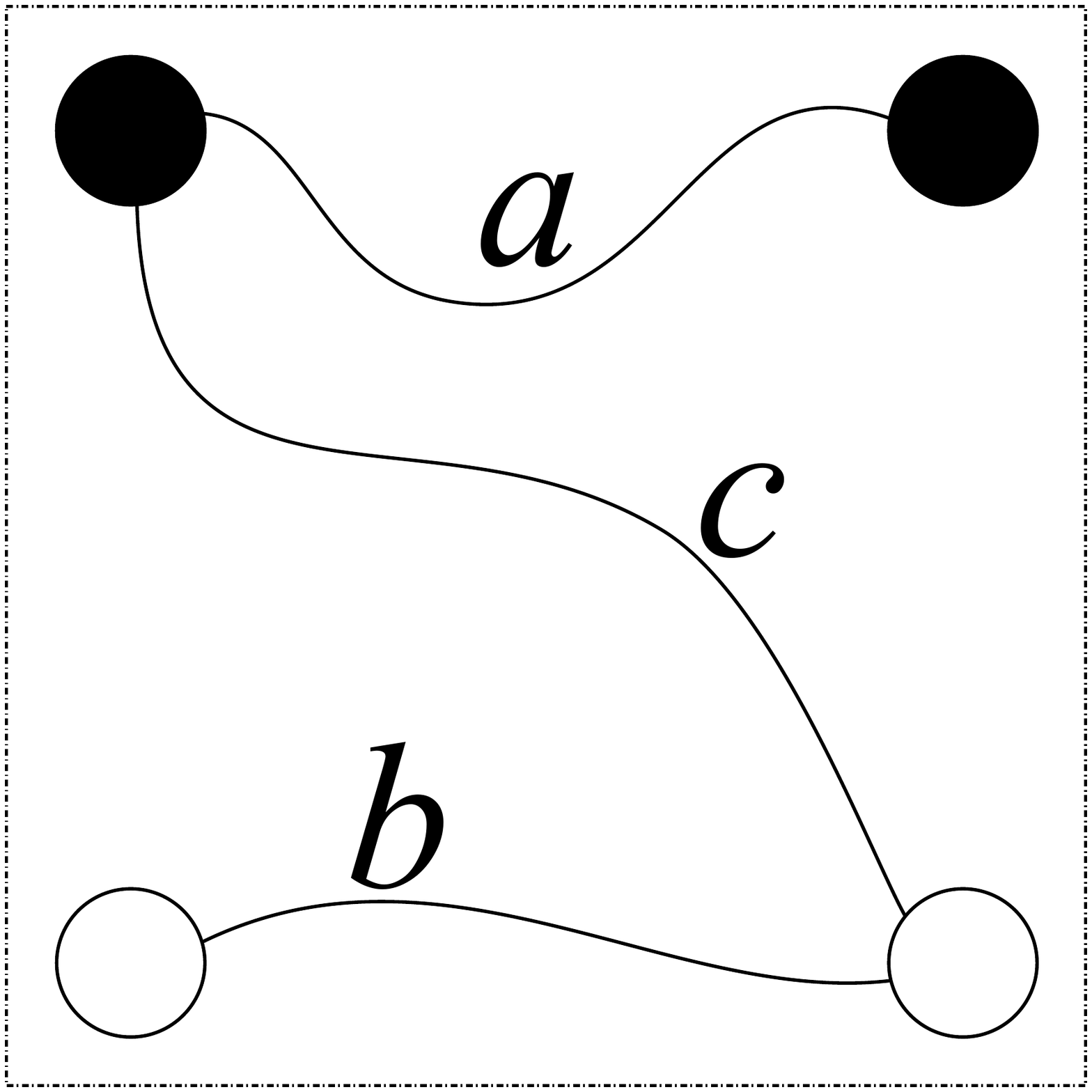}
\centering
\caption{Assumption of Lemma 7.}
\centering
\end{figure}

Call the subchains between the corner nodes branches $a$, $b$, and $c$. Branches $a$ and $b$ must be all G and C, respectively. If either of those branches does not follow straight against the wall of the bounding box, unbound nodes will be created and the folding will not be optimal.

\begin{figure}[h]
\includegraphics[scale=0.15]{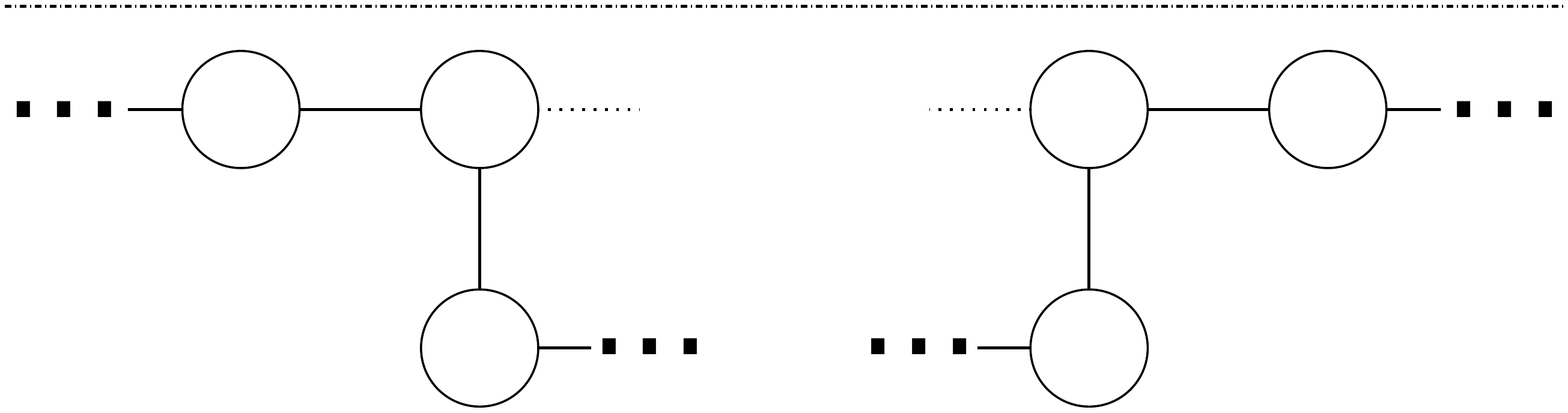}
\centering
\caption{Branch $a$ if it does not follow the wall.}
\centering
\end{figure}

The bounding box is rectangular and $a$ and $b$ are straight lines, so $c$ is symmetrical (half G and half C). Examine the corner where $b$ and $c$ meet.

\begin{figure}[h]
\includegraphics[scale=0.15]{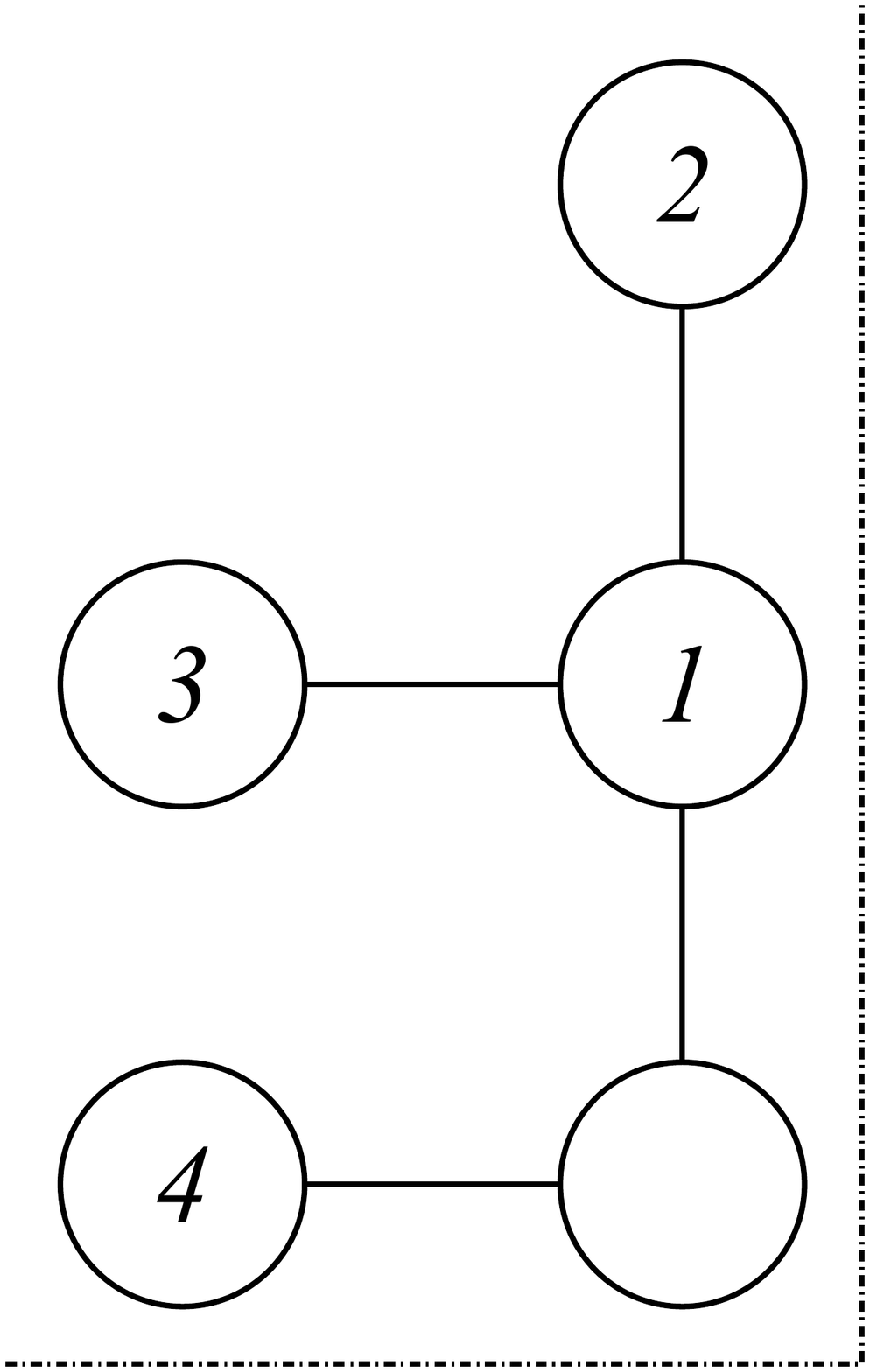}
\centering
\caption{Illustration of corner between $a$ and $b$.}
\centering
\end{figure}

After following branch $b$ to the corner of the bounding box, node $1$ is forced to be in its position, and there are two choices where to put the next node (which must also be a C because $n>5$), either $2$ or $3$. If the node is put in position $3$, node $4$ is unable to bond, and the folding is not optimal. If the next node is put in position $2$ and a G node fills position $3$, both nodes $1$ and $4$ need to bond with it, but only one of them can so the folding is not optimal. Therefore the assumption is wrong.  $\Box$
\newline \indent
Note: this proof only holds as long as the two nodes after the corner are the same assignment as the corner, but that is true whenever $n>3$.

\bigskip
\noindent
\textbf{Proof of Theorem 1.} By Fact 2 and Lemma 7, an optimal folding of $S_n$ must resemble Figure 8.

\begin{figure}[h]
\includegraphics[scale=0.13]{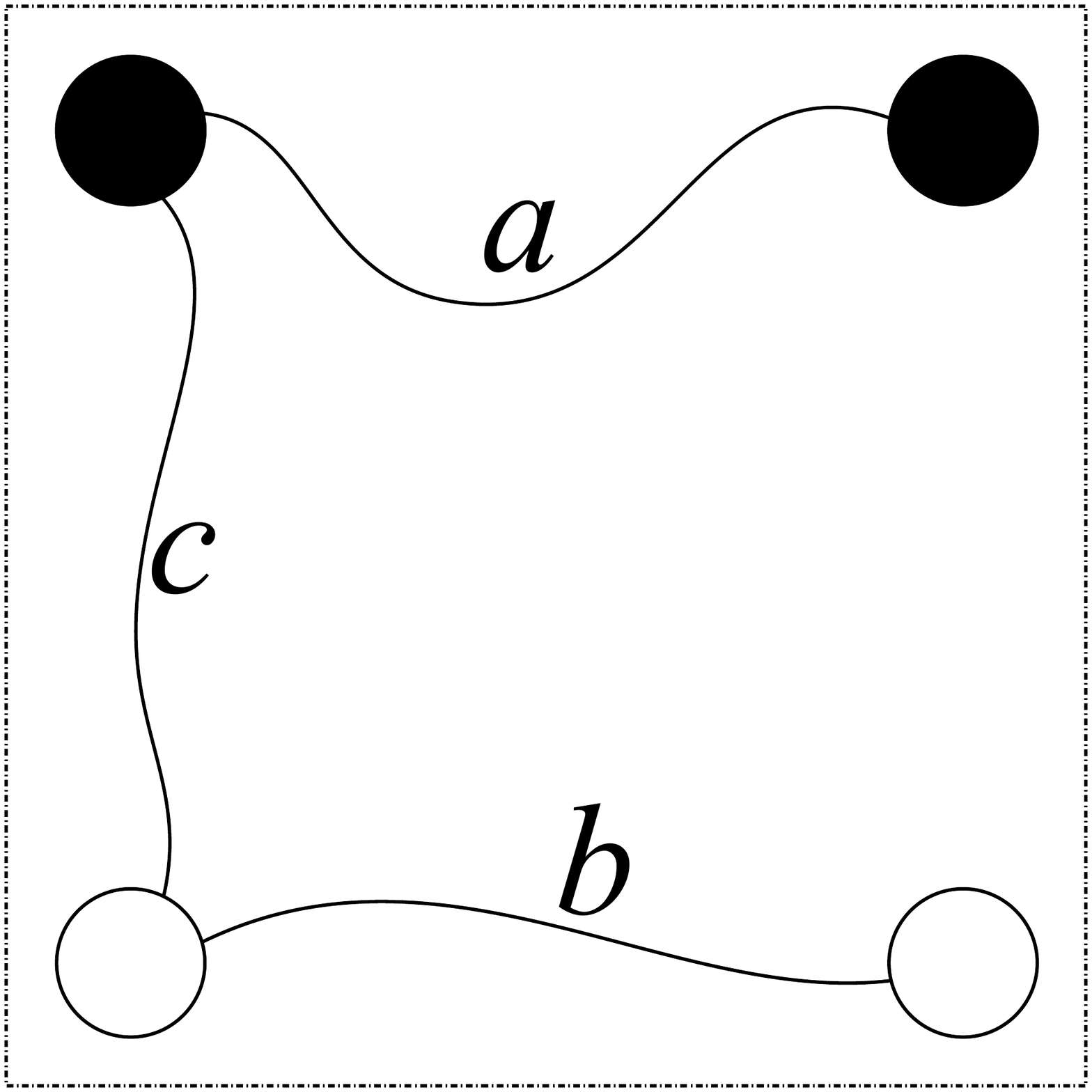}
\centering
\caption{Approximate folding of $S_n$.}
\centering
\end{figure}

Just as in the proof of Lemma 7, branches $a$ and $b$ must follow the wall of the bounding box and $c$ is symmetrical. Assume $|c| \neq 0$. The nodes of $c$ must bond to every node of both $a$ and $b$, including the endpoints. However, all of the G nodes on $c$ are on the top and all the C's are on the bottom. Therefore, as soon as a single G node on $c$ bonds to $b$, all of the C nodes on $c$ are trapped and cannot bond to $a$. This means that $|c|=0$ and $F_n$ is the uniquely optimal folding of $S_n$. $\Box$

\bigskip
\noindent
\textbf{Corollary 8.} Theorem 1 assumes that there are no A or U bases in the sequence. In reality all four bases are used, so it would be preferable to be able to show that a more even mix of bases could still have a unique optimal folding. If the desired ratio of $G/C:A/U = m:n$, consider the chain $S_n=G^{m/2}A^{n/2}U^{n/2}C^{m/2}$. Temporarily allow G-U and A-C bonds. Apply Theorem 1 to show that a unique optimal folding exists under these conditions. Reinstating Watson-Crick pairing does not break any bonds of the folding, and it can only cut down on the number of bonds in other foldings. Therefore $S_n$ must still have a unique optimal folding. $\Box$

\section{NP Hardness}
\label{sec4}

\noindent
\textbf{Theorem 2.} \textit{For a chain $F_n$ and a number k, finding a folding of $F_n$ with at least k bonds is NP-hard.}

\bigskip
\noindent
Finding a suitable folding will be reduced from rectilinear planar monotone 3-SAT by the same method as Theorem 2 of \cite{flattening}, which describes a polynomial-time algorithm for reducing from rectilinear planar monotone 3-SAT to flattening a fixed angle chain. To do this, only three gadgets are needed: a forced direction 90\degree turn, a 90\degree turn where the direction can be chosen, and a straight line. These are all attainable using only 5 unique bases: two pairs and one base which can not bond. The entire linkage will be represented by a double strand of RNA. Straight portions that can flip directions have the pattern $(CCCA)^n/(GGGU)^n$ (4-cycle), and ones that cannot have the pattern $(CCCCCCCA)^n/(GGGGGGGU)^n$ (8-cycle). Fixed direction edges are shown with arrows. Note that when turning left or right, the alignment of the two strands becomes shifted by two bases relative to the straight portion, and that there are two unbound nodes at every variable turn. Therefore, $k=n/2-t$ where $n$ is the number of bondable nodes and $t$ is the number of variable turns. 

\begin{figure}[h]
\includegraphics[scale=0.17]{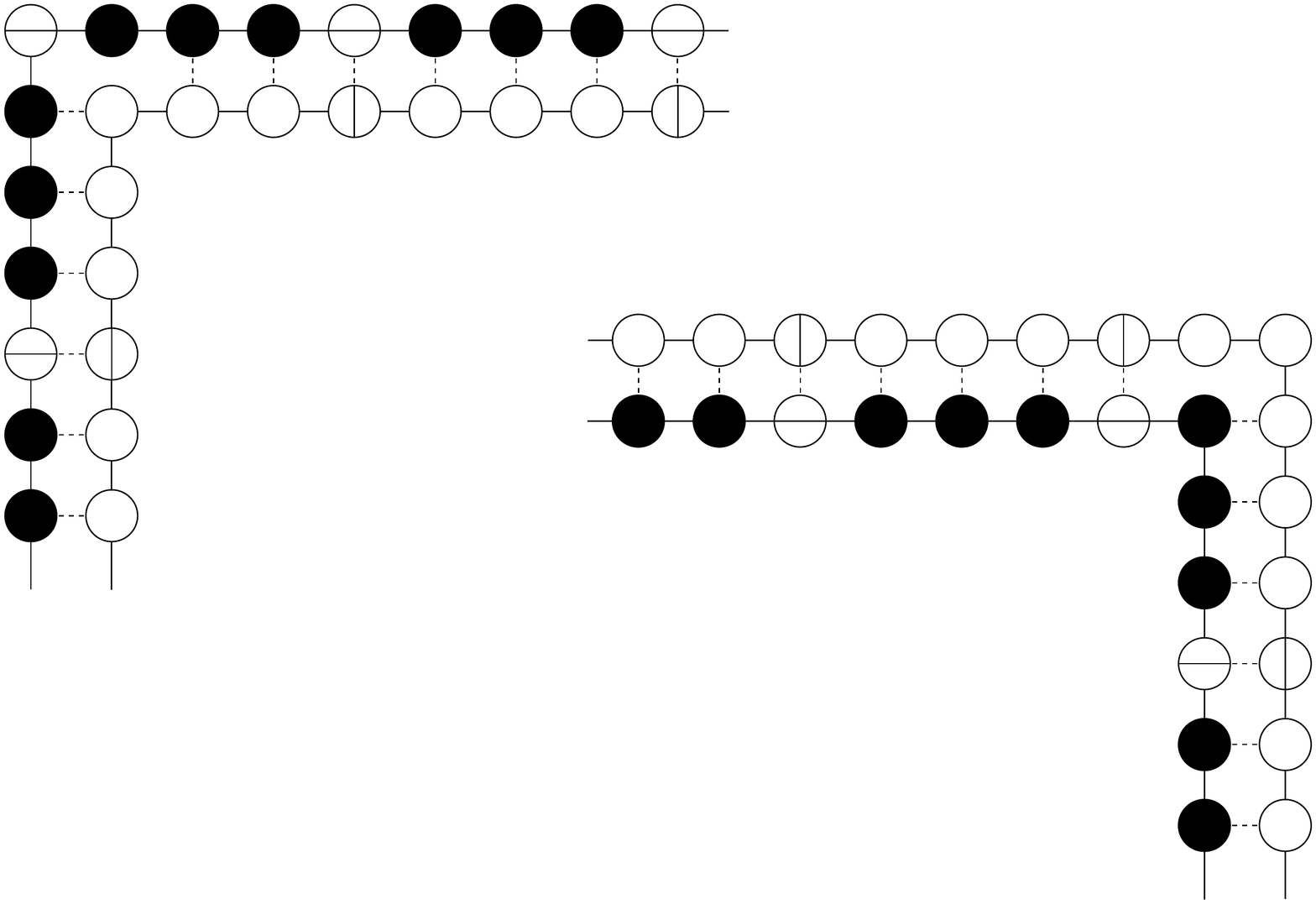}
\centering
\caption{Variable turn gadget, represented by a dot on a vertex.}
\centering
\end{figure}

\begin{figure}[h]
\includegraphics[scale=0.17]{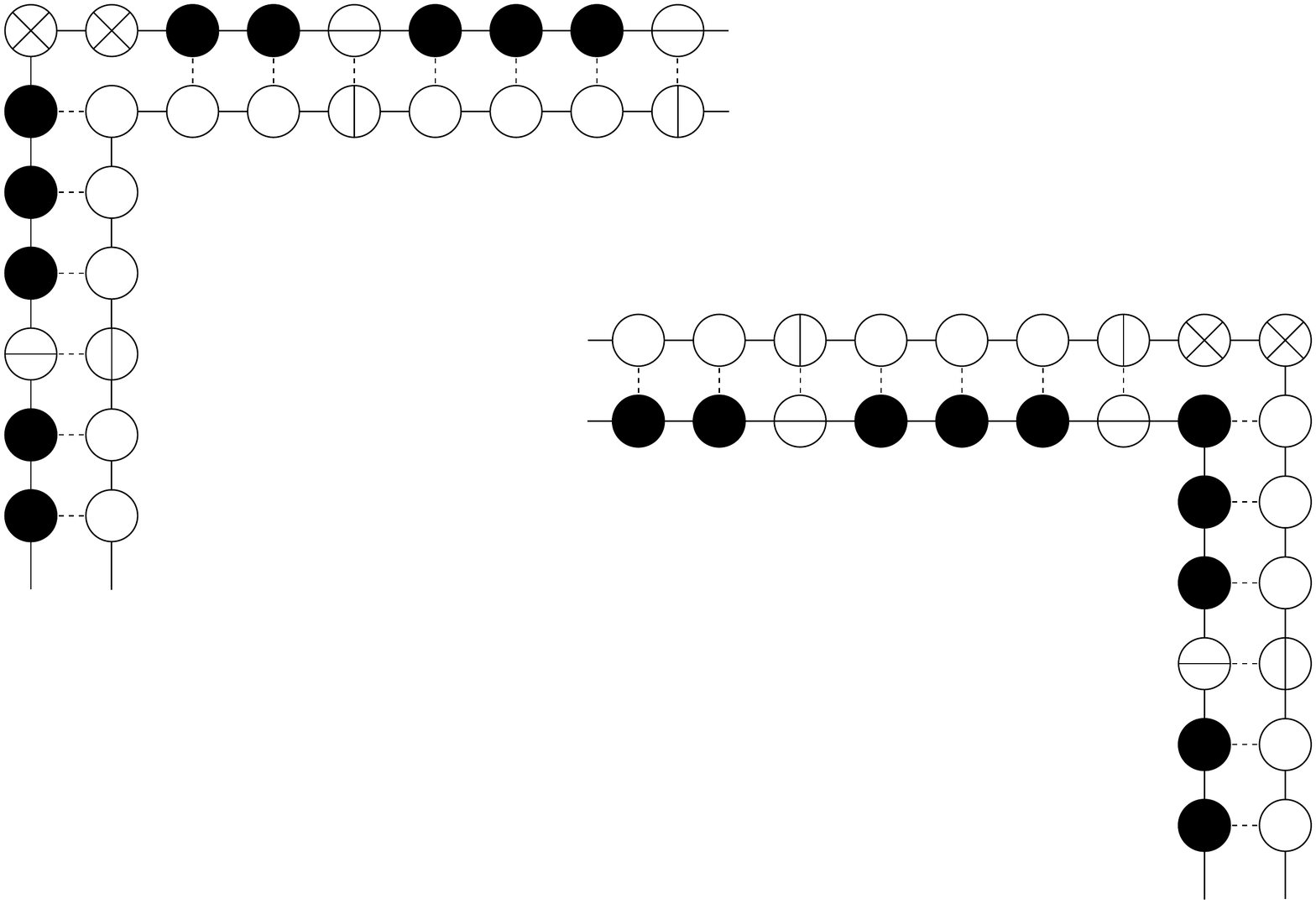}
\centering
\caption{Fixed turn gadgets, which are shown with no dot.}
\centering
\end{figure}

Using these gadgets, we will create the same clause and variable structure as in [Demaine, Eisenstat 2011]. There is one major modification that needs to be made. In the 2011 paper, a variable turn causes the subsequent chain to reflect over the preceding edge, while in the RNA model a turn causes the rest of the chain to rotate by 180\degree. However, note that every variable turn in [Demaine, Eisenstat 2011] has a complementary turn which reverses the reflection caused by the first. There are occasionally constructions between the turns, but those constructions do not interfere with the pair of turns.

\begin{figure}[h]
\includegraphics[scale=0.17]{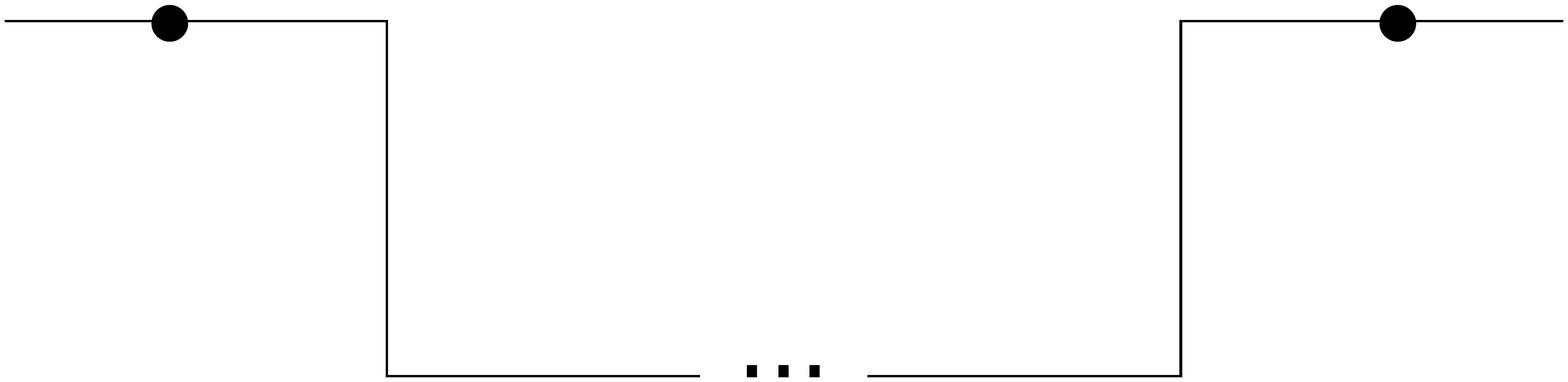}
\centering
\caption{Turns are always paired.}
\centering
\end{figure}

To replicate this behavior with RNA, use variable turns and fixed direction edges:

\begin{figure}[h]
\includegraphics[scale=0.17]{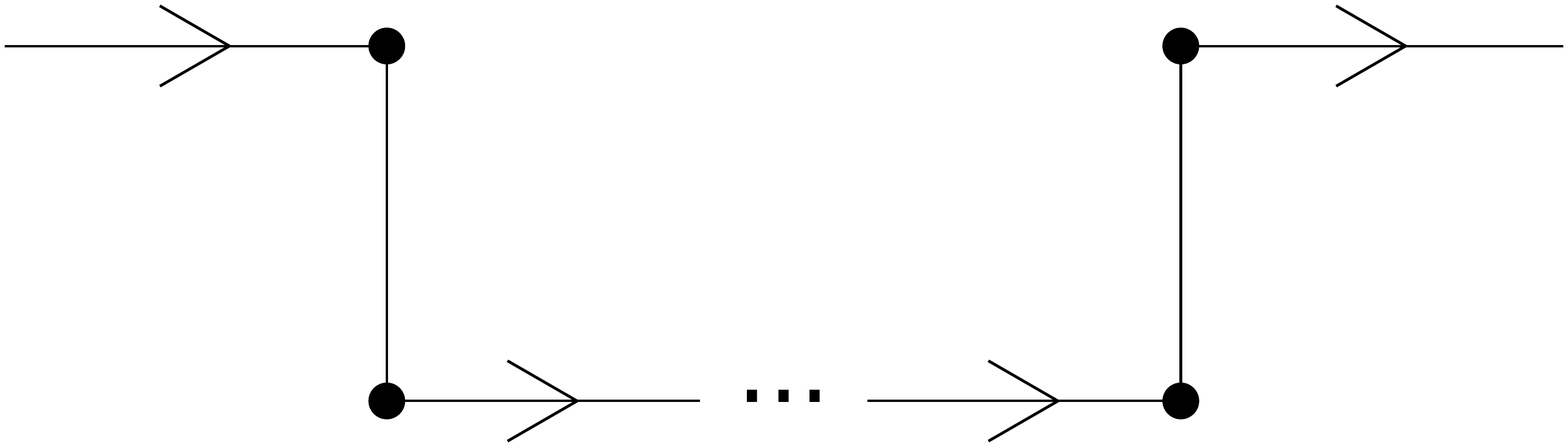}
\centering
\caption{RNA version of a pair of turns.}
\centering
\end{figure}

This is enough to replicate the structure and intended function of the gadgets in \textit{Flattening Fixed-Angle Chains
Is Strongly NP-Hard}, but it still needs to be shown that unintended foldings are not possible.

First, a few notes on constructing the reduction. Because the alignment of the strands can be shifted by 4 at every variable turn, it is a concern that eventually the chains could intersect even if a valid 3-SAT was found. To avoid this, ensure that the spacing between each chain is greater than $40t$, where $t$ is the number of variable turns. This only increases the size of the chain by a polynomial factor, so the reduction from 3-SAT still holds. By the same logic, is does not matter that each turn can be made at three positions around the intended point without losing any bonds. In proportion to the entire structure, such a small change is insignificant.

The variable portion of the construction will also be doubled up so that the two endpoints meet up. The overview of the structure is shown in Figure 14.

Two long double stranded tails of X nodes are attached to the endpoints, pointing away from the rest of the folding. These segments each have at least $(N/2)^2$ nodes in them, where $N$ is the total number of nodes in the rest of the chain.

\begin{figure}[h]
\includegraphics[scale=0.17]{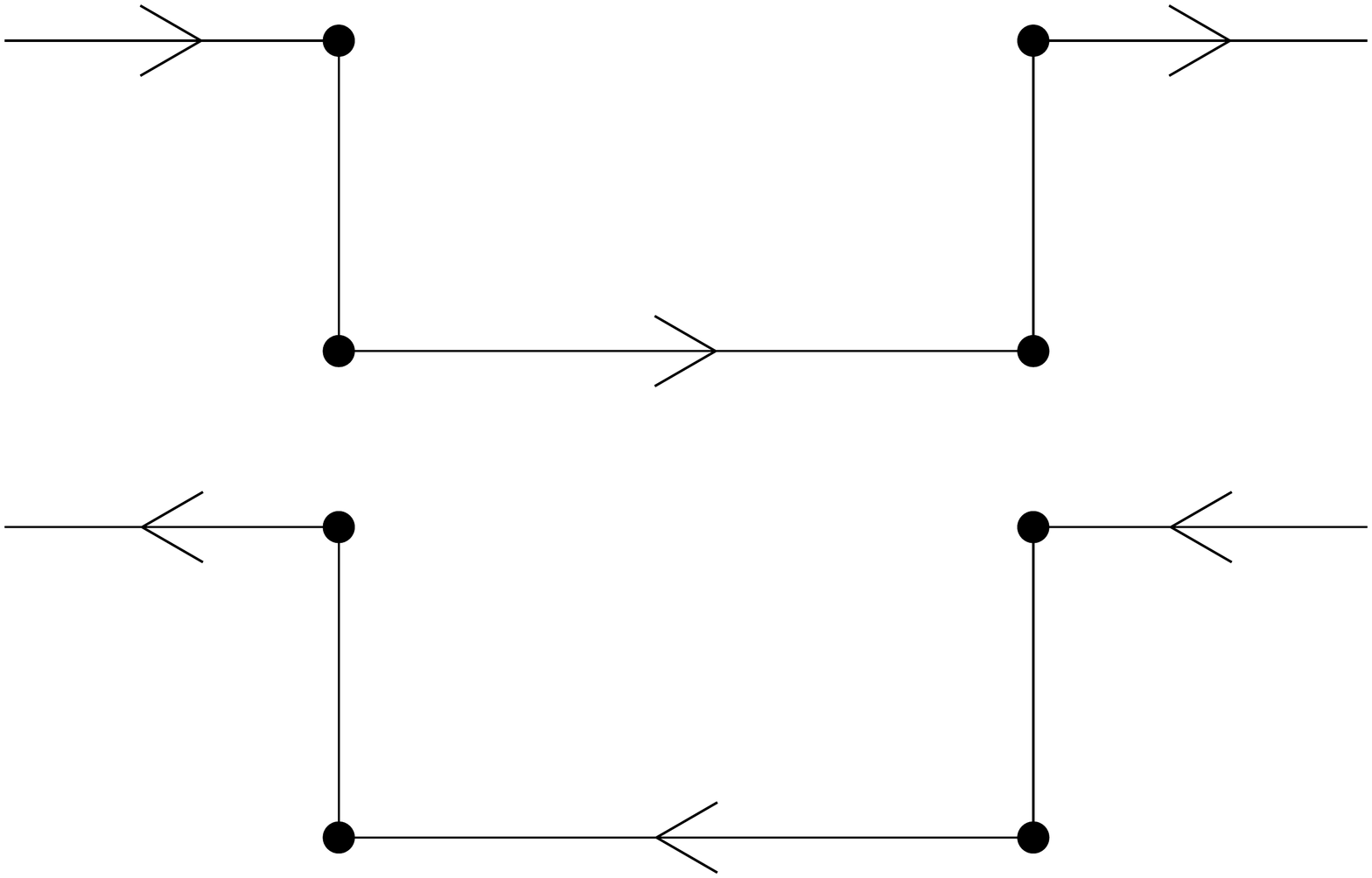}
\centering
\caption{A variable gadget.}
\centering
\end{figure}

\begin{figure}[h]
\includegraphics[scale=0.17]{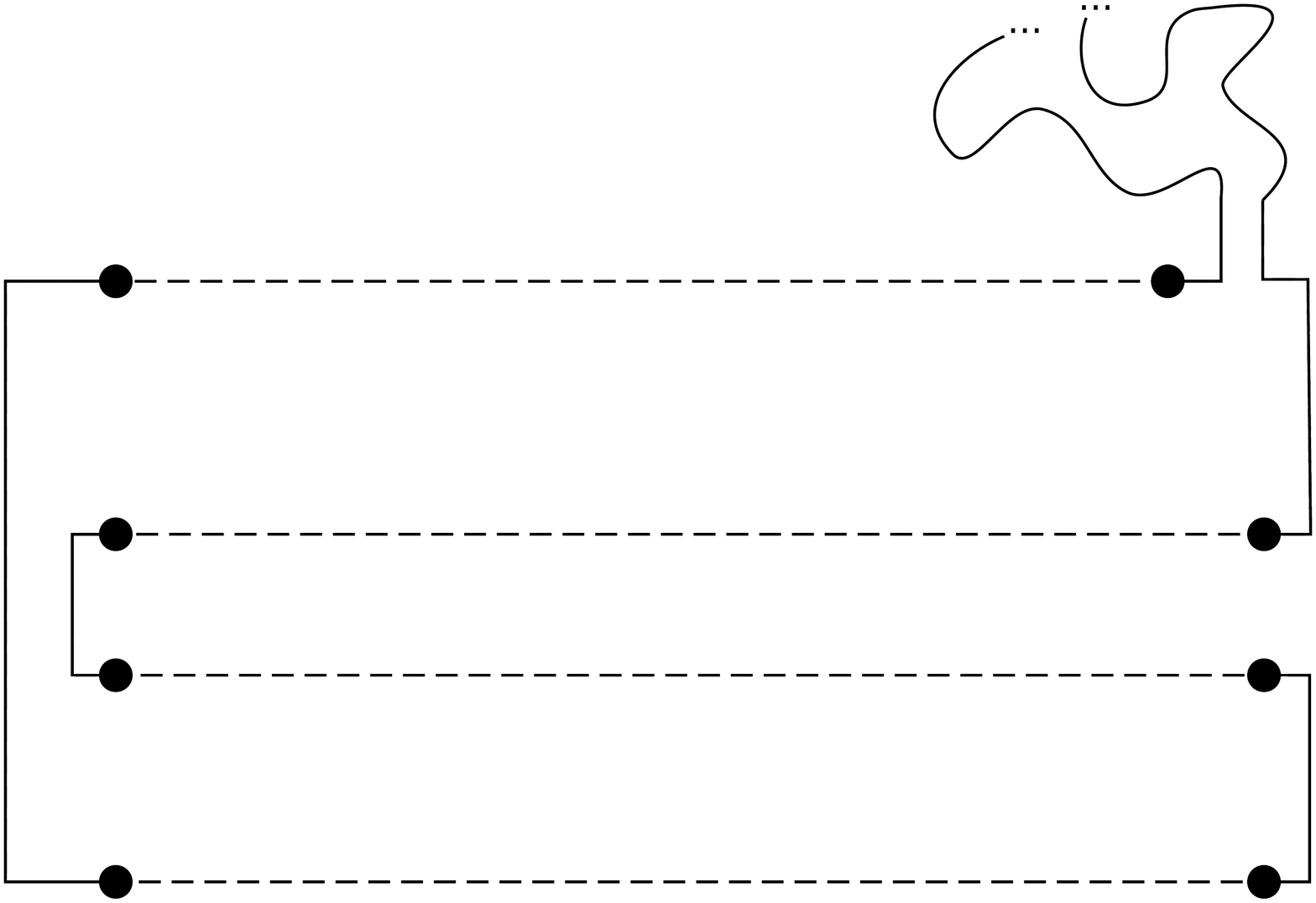}
\centering
\caption{Large scale structure of the construction, with the chain of the variable gadgets between the two chains of clause gadgets, and the X tails at the top.}
\centering
\end{figure}

\begin{figure}[h]
\includegraphics[scale=0.2]{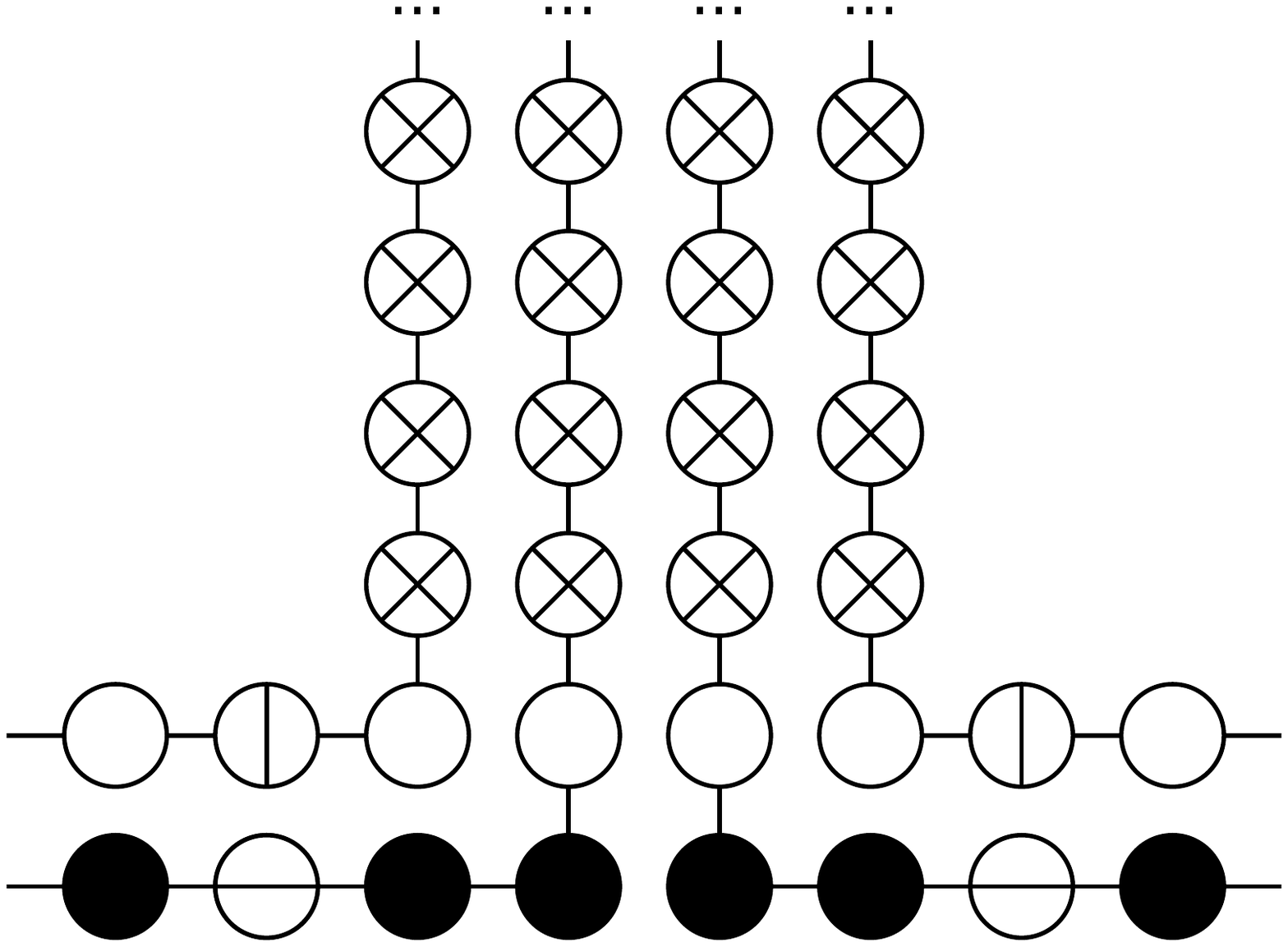}
\centering
\caption{Endpoints with chains of X nodes}
\centering
\end{figure}

\bigskip
\noindent
\textbf{Lemma 9.} \textit{Isolated straight sections must remain straight.}

\bigskip
\noindent
\textbf{Proof.} Apply Theorem 1 to show that a straight double strand is the only optimal folding of $(CCCA)^n/(GGGU)^n$ or $(CCCCCCCA)^n/(GGGGGGGU)^n$. Any other folding will lose bonds, and so in order to still achieve $k$ total bonds must form other bonds somewhere else. However the section is isolated so this is impossible. $\Box$

\bigskip
\noindent
\textbf{Lemma 10.} \textit{No unintended foldings are possible.}

\bigskip
\noindent
\textbf{Proof.} All of the A and C nodes are on one side of the chain, and U and G on the other. Therefore, locally the only bonds that can form are across the double strand. On a larger scale, more bonds could be formed if the chain wrapped around the endpoint, bringing the A/C chain into contact with the U/G chain. However, the X tails at each endpoint each have at least $(N/2)^2$ nodes, which means the chain can not wrap around them and bond with the unbound nodes left on variable turns. Any fold or alteration that is not made at a designated corner will lose bonds, and because the goal number of bonds only takes into account those made at intended corners (it is impossible not for lose bonds at intended corners), only intended foldings can possibly achieve that number. $\Box$

\bigskip
\noindent
\textbf{Proof of Theorem 2.} Constructing a RNA sequence from an instance of rectilinear planar monotone 3-SAT in polynomial time is possible, and the only way to find a folding with at least $k$ bonds is by solving the 3-SAT, so finding finding a folding of a given RNA sequence with at least $k$ bonds is NP-hard. $\Box$

\bigskip
\noindent
\textbf{Notes.} This construction uses five bases: A, G, C, U, and X. Only four of these form bonds. Using a constant number of types of bases is highly desirable because it emulates nature as closely as possible.

\section{Constant-Factor Approximation Algorithm}
\label{sec5}

\noindent
\textbf{Theorem 3.} \textit{A linear time algorithm exists which will find a folding of a given RNA sequence made up of only C and G bases with at least 1/12 of the optimal number of bonds.}

\bigskip
\noindent
This algorithm is a simple modification to the one used to achieve a 1/3 approximation of folding proteins in the HP model done in \textit{A New Algorithm for Protein Folding in the HP Model}\cite{approximation}.

\textbf{Proof of Theorem 3.} Although there are only two node assignments, there are actually four types of nodes because of the parity that occurs on a square grid. C and G nodes can only bond if they have opposite parity. We call these odd and even nodes based on their order in the original chain. An upper bound on the total number of bonds that can be formed is $min(odd G, even C)+min(even G, odd C)$ (where $odd G$ is the number of odd-G's) by a simple parity argument. To map these four types onto the algorithm of [Newman 2002], we do the following: If $min(odd-G, even C)\geq min(even G, odd C)$, we rename odd-G's odd-1's, even-C's even-1's, and even-G's and odd-C's 0's. Otherwise, we rename even-G's even-1's, odd-C's odd-1's, and odd-G's and even-C's 0's. Note that this loses at most 1/2 of the optimal number of bonds. Next, we apply the HP model approximation algorithm, which guarantees 1/3 of the maximum\cite{approximation}. However, we again lose another 1/2 of available bonds because the HP model assumes that each node can form up to two bonds whereas the RNA model restricts that number to 1. Multiplying these factors together shows that this method will form at least 1/12 of the optimal number of bonds for a given RNA sequence. $\Box$

\section{Acknowledgements}
I thank Erik Demaine at the MIT Computer Science and Artificial Intelligence Laboratory for his extensive support and mentoring during the course of this research, and Julie Blackwood at Williams College for her guidance in editing and submitting this paper.

%\appendix

%\section{Section in Appendix}
%\label{appendix-sec1}

%Sample text. Sample text. Sample text. Sample text. Sample text. Sample text. 
%Sample text. Sample text. Sample text. Sample text. Sample text. Sample text. 
%Sample text. 

%% References


\bibliographystyle{elsarticle-num}

\bibliography{RNAFolding}

\end{document}